\documentclass[a4paper,12pt]{article}
\usepackage{amsfonts}
\usepackage{amsthm}
\usepackage{amsmath} 
\usepackage{amsfonts}
\usepackage{latexsym}
\usepackage{amssymb}

\newtheorem{teo}{Theorem}[subsection]
\newtheorem{lem}[teo]{Lemma}
\newtheorem{cor}[teo]{Corollary}
\newtheorem{pro}[teo]{Proposition}
\newtheorem{fed}[teo]{Definition}
\newtheorem{rem}[teo]{Remark}
\newtheorem{conj}{Conjecture}

\newtheorem{claim}[teo]{Claim}

\def\noi{\noindent}
\def\bdem{\begin{proof}}
\def\edem{\renewcommand{\qed}{\hfill $\blacksquare$}\end{proof}}
\def\QED{\hfill $\blacksquare$}
\def\EOE{\hfill $\blacktriangle$}

\def\eps{\varepsilon}
\def\fii{\varphi }
\def\la{\lambda}
\def\cF{\mathcal{F}}

\def\al{\alpha}

\def\si{\sigma}

\def\ds{\displaystyle}

\def\N{\mathbb{N}}

\def\R{\mathbb{R}}
\def\C{\mathbb{C}}
\def\inc{\subseteq}
\def\bm{\left(\begin{array}}
\def\em{\end{array}\right)}

\def\cA{\mathcal{A}}

\def\cE{\mathcal{E}}

\def\cO{\mathcal{O}}
\def\cP{\mathbb{P}}

\def\cM{\mathcal{M}}
\def\cN{\mathcal{N}}

\def\cU{\mathcal{U}}
\def\cV{\mathcal{V}}
\def\cW{\mathcal{W}}

\def\ese{\mathcal{S}}
\def\ete{\mathcal{T}}
\def\ewe{\mathcal{W}}

\def\D{\Delta}

\def\ben{\begin{enumerate}}
\def\een{\end{enumerate}}
\def\beq{\begin{equation}}
\def\eeq{\end{equation}}
\def\barr{\begin{array}}
\def\earr{\end{array}}
\def\inv{^{-1}}
\def\sii{ if and only if }

\newcommand{\peso}[1]{ \quad \text{ #1 } \quad }
\newcommand{\sub}[2]{{#1}_{\mbox{\tiny{${#2}$}}}}

\DeclareMathOperator{\Preal}{\R\mbox{e}} 

\DeclareMathOperator*{\dist}{dist}
\DeclareMathOperator*{\convsotdpre}{\searrow}

\DeclareMathOperator{\tr}{tr}

\DeclareMathOperator{\leqp}{\leqslant}

\def\al{\alpha}

\newcommand{\pint}[1]{\displaystyle \left \langle #1 \right\rangle}

\newcommand{\hil}{\mathcal{H}}

\newcommand{\cene}{\mathbb{C}^r}
\newcommand{\mat}{\mathcal{M}_r (\C) }
\newcommand{\matu}{\mathcal{U}(r)}
\newcommand{\matsa}{\mathcal{M}_r^{h}(\C)  }
\newcommand{\matah}{\mathcal{M}_r^{ah}(\C)  }

\newcommand{\matinv}{\mathcal{G}\textit{l}\,_r(\C) }

\newcommand{\spec}[1]{\sigma\left( #1\right)}

\newcommand{\conv}{\xrightarrow[n\rightarrow\infty]{}}

\newcommand{\convd}{\convsotdpre_{n\rightarrow\infty}^{}}

\newcommand{\alu}[1]{\Delta\left(#1\right)}
\newcommand{\alul}[1]{\Delta_\la\left(#1\right)}
\newcommand{\aluf}[1]{\left|#1\right|^{1/2}U\left|#1\right|^{1/2}}
\newcommand{\aluit}[2]{\Delta^{#1}\left(#2\right)}
\newcommand{\alulit}[2]{\Delta_\la^{#1}\left(#2\right)}

\newcommand{\orb}[1]{\ese \left({#1}\right)}
\newcommand{\orbu}[1]{\cU \left({#1}\right)}
\newcommand{\der}[3]{\sub{T}{#2}{#1} \left({#3}\right)}
\newcommand{\dersin}[2]{\sub{T}{#2}{#1}}

\newcommand{\mla}[1]{\sub{\cM}{{\la,{#1}}}}

\newcommand{\lentas}{\cO_D}

\newcommand{\cubo}[1]{\mathcal{Q}_{#1}}

\begin{document}
\title{\textbf{The iterated Aluthge transforms of a matrix converge}
\author{Jorge Antezana\thanks{Partially supported by CONICET (PIP 4463/96), Universidad de La Plata (UNLP 11 X472) and
ANPCYT (PICT03-09521).}\and Enrique R. Pujals \thanks {Partially supported by CNPq} \and Demetrio Stojanoff$^{\,*}$ }
\date{}
}
\maketitle

\vglue1.5truecm

\begin{abstract}
Given an $r\times r$ complex matrix $T$, if $T=U|T|$ is the polar
decomposition of $T$, then, the Aluthge transform is defined by 
$$
\Delta\left(T \right)= |T|^{1/2} U |T |^{1/2}.
$$ 
Let $\Delta^{n}(T)$ denote the n-times iterated Aluthge transform of $T$, i.e.
$\Delta^{0}(T)=T$ and $\Delta^{n}(T)=\Delta(\Delta^{n-1}(T))$, $n\in\mathbb{N}$. 
We prove that the sequence $\{\Delta^{n}(T)\}_{n\in\mathbb{N}}$ converges 
for every $r\times r$  matrix $T$. This result was conjecturated by 
Jung, Ko and Pearcy in 2003. 
We also analyze the regularity of the limit function. 
\end{abstract}

\vglue1truecm

\noi

\medskip
\noi

\section{Introduction}

Let $\hil$ be a Hilbert space and $T$ a bounded operator defined on $\hil$ whose polar 
decomposition is $T=U|T|$. The \textit{Aluthge transform} of $T$ is the operator 
$\alu{T}=|T|^{1/2}U\ |T|^{1/2}$. This transform was introduced in \cite{[Aluthge]} 
to study p-hyponormal and log-hyponormal operators. Roughly speaking, the idea 
behind the Aluthge transform is to convert an operator into other operator which shares 
with the first one some spectral properties but it is closer to being a normal operator.

The Aluthge transform has received much attention in recent years. One reason is its connection 
with the invariant subspace problem. Jung, Ko and Pearcy proved in \cite{[JKP0]} 
that $T$ has a nontrivial invariant subspace if an only if $\alu{T}$ does. On the 
other hand, Dykema and Schultz proved in \cite{[Dykema]} that the Brown measure 
is preserved by the Aluthge transform.

Another reason is related with the iterated Aluthge transform. 
Let $\aluit{0}{T}=T$ and $\aluit{n}{T}=\alu{\aluit{n-1}{T}}$ 
for every $n\in\N$. In  \cite{[JKP1]} Jung, Ko and Peacy raised the following conjecture:

\begin{conj}
\label{JKP} \rm
The sequence  of iterates $\{\aluit{n}{T}\}_{n\in\N}$ 
 converges, for every 
 matrix $T$.
\end{conj}

\noi Although many results supported this conjecture (see for instance \cite{[Ando]} and 
\cite{[Yamazaki]}), there were only partial solutions. For instance, 
Ando and Yamazaki proved in \cite{[Ando-Yamaza]} that Conjecture \ref{JKP}
is true for $2\times 2$ matrices, Dykema and Schultz in \cite{[Dykema]} proved 
that the conjecture is true for an operator $T$ in a finite factor such that 
the unitary part of its polar decomposition normalizes an abelian subalgebra 
that contains $|T|$, and Huang and Tam proved in \cite{[HT]} that the conjecture is true for
matrices whose eigenvalues have different moduli. 

In our previous work \cite{[APS]}, we introduced a new approach to studying this problem, which was based 
on techniques from dynamical systems. This approach allowed to show that 
Conjecture \ref{JKP} is true for every diagonalizable matrix. 

In this paper, using again 
dynamical techniques, combined with some geometrical arguments, we completely solve 
Conjecture \ref{JKP}. 
There are many fruitful points of contact between the theories of dynamical systems and operator algebras. The combination of dynamical and geometrical techniques used to study Conjecture \ref{JKP} suggests a new possible interaction between both theories. On one hand, it provides another field of applications of the stability theory of hyperbolic systems and invariant manifolds. In this sense, the work by Shub and Vasquez on the $QR$ algorithm is an important 
precedent (see \cite{Sh3}). On the other hand, it provides to the operator theorists a new set of powerful tools to deal with problems where the usual techniques fail. In our case, besides the solution, it also provides a better understanding of  the problem. The dynamical 
perspective not only 
allows to prove Conjecture \ref{JKP}, but it also provides further information related to the
regularity of the limit function and the rate of convergence of the iterated sequence.

By a result proved in \cite{[AMS]}, Conjecture \ref{JKP} reduces to the invertible case. 
On the other hand, according to a result independently 
proved by Jung, Ko and Pearcy in \cite{[JKP1]}, and  by Ando in \cite{[Ando]}, 
for every  invertible $r\times r$ matrix 
$T$, the sequence of iterates $\{\aluit{n}{T}\}$ goes toward  the set of normal operators 
which have the same characteristic polynomial as $T$. This set can be characterized 
as the unitary orbit of some matrix $D$ that also shares the characteristic 
polynomial with $T$. Let $\orbu{D}$ denote this unitary orbit. 

If $T$ is diagonalizable, then $T$ and all the iterates $\aluit{n}{T}$ 
belong to the similarity orbit of $D$, denoted by $\orb{D}$, 
which is a riemannian manifold that contains $\orbu{D}$ 
as a compact submanifold. 
Note that all the points of $\orbu{D}$ are normal matrices, and therefore they 
are fixed points for the Aluthge transform. These facts suggest the possibility 
of pursuing a dynamic approach in $\orb{D}$ based in the stable manifold theorem. 
This approach was carried out in  \cite{[APS]} to prove the convergence 
of the iterated Aluthge transform sequence for diagonalizable matrices. 

However, the non-diagonalizable case is 
different, since 
the geometry context of the problem is more complicated. Indeed, 
if $T$ is not diagonalizable, $\orbu{D}$ is contained in the boundary of $\orb{T}$, 
which also contains the orbits of matrices with smaller Jordan forms than 
the Jordan form of 
$T$. The boundary of $\orb{T}$ can be thought as a sort 
of lattice of boundaries. Therefore, in order to prove Conjecture \ref{JKP}, 
we put the problem in a different setting so that both cases can be 
analyzed together. 
The Aluthge transform is viewed as an endomorphism on the space of invertible 
matrices,  and we consider all the orbits mentioned before
not as a manifold, but as the basin of attraction $B_\Delta(\orbu{D}\,)$ i.e.,  
those matrices $T$ such that the sequence $\{\aluit{n}{T}\}_{n\in\N}$ goes to 
$\orbu{D}$  as $n\conv \infty$. 
The basin $B_\Delta(\orbu{D}\,)$ can also be characterized 
as the set of those matrices that have the same characteristic polynomial as $D$.

The stable manifold theorem can be 
extended to $B_\Delta(\orbu{D}\,)$ (see Thm. \ref{pelitosA}),  and 
no differential structure is required in the basin. Using this
theorem we construct, through each 
$T$ in the basin close enough to $\orbu{D}$, a $\Delta$-invariant manifold $\ewe^{ss}_T$ 
which satisfies that 
$$
\ewe^{ss}_T\subseteq \{S: \|\aluit{n}{T}-\aluit{n}{S}\|< C\gamma^n \ 
\mbox{ for every } \ n \in \N \}\,,
$$
where $C$ and $\gamma<1$ are constants that only depend on the distance among different eigenvalues of $D$. 
Hence, if the sequence $\aluit{\infty}{S}$ converges for  some  $S \in \ewe^{ss}_T \,$, then the same must 
happen for $T$.
For the diagonalizable case, in \cite{[APS]} we have considered only the stable 
manifolds $\ewe^{ss}_N$ for points $N \in \orbu{D}$. Then, 
using an argument which involves the inverse mapping theorem, we deduced that the union of these 
manifolds contains an open neighborhood of $\orbu{D}$ in $\orb{D}$. 

That approach fails in the general case, because the basin in not a manifold. In this case we prove that, 
for every $T$ in the basin near to $\orbu{D}$, the stable manifolds $\ewe^{ss}_T$  intersect the set  
$\cO_D\,$ of those matrices in the basin with orthogonal spectral projections.  The set $\cO_D\,$ 
 can be studied with the usual properties of the Aluthge transform, showing that
Conjecture \ref{JKP} holds for its elements. 
However,  $\cO_D$ does not have a differential structure. 
To avoid this problem, in order to see that $\ewe^{ss}_T \cap \cO_D\neq \varnothing $,  
we project the stable  manifold $\ewe^{ss}_T$ and $\cO_D$ 
to the orbit $\orb{D}$, using spectral projections, in such a way that $\cO_D$ projects 
onto the manifold $\orbu{D}$. 
Then, inside the manifold $\orb{D}$, the desired result follows by a geometric argument 
based in known results about transversal intersections.  

Another problem arises  proving 
the continuity of the limit map 
$\Delta^{\infty}$ on $\matinv$:
If $N_0 \in \matinv$ is normal, but has eigenvectors with multiplicity greater that 1, then 
in every open neighborhood of $N_0$ there exist matrices with very close but different eigenvalues.
This fact reduces drastically the rate of convergence for such matrices, 
even in the diagonalizable case (see section 6.1). To solve this problem we separate  
the spectrum of these matrices in blocks which are near to each eigenvalue of $N_0\,$, 
even if in these blocks the eigenvalues are different. Then, we repeat the strategy 
of the proof of the convergence, but with respect to spectral projections relative 
the blocks indexed by the spectrum of $N_0$ (instead of using the spectrum of each 
matrix near $N_0\,$). As before, we show that these projections converge 
at an uniform velocity to an orthogonal system of projections, near the system of
$N_0\,$. Then, one can see easily that,  whatever were the limit and the rate of convergence
to this limit, it must remain close to $N_0\,$, because the (block) spectral projections  
and the global spectra are close.

Using the results of \cite{[APS2]}, all the results mentioned before 
can be extended to the so called $\la$-Alutghe transform 
$\alul{T} = |T|^{\la}U\ |T|^{1-\la}$, for every $\la \in (0,1)$.

This paper is organized as follows: 
in section 2, we collect several preliminary definitions and results about 
the  Aluthge transform, the geometry of similarity 
and unitary orbits, the stable manifold theorem in local basins,  
and the known properties of the spectral projections. 
In section 3, we compute the derivative of $\Delta$ in the whole space 
$\mat$, and we state the Dynamical Systems aspects of $\Delta$. Particularly 
the stable manifold theorem on the local basin of a compact set of 
fixed points. In section 4 we prove Conjecture \ref{JKP} about the convergence 
$\aluit{n}{T} \conv \aluit{\infty}{T}$ for every $T \in \mat$. In section 
5 we study the regularity of the limit map $T \mapsto \aluit{\infty}{T}$, mainly for 
$T$ invertible. Section 6 contains concluding remarks about the rate of convergence, 
 and the extension of the main results 
to the $\la$-Aluthge transforms, for every $\la \in (0,1)$. 
In the Appendix, we write the proof of two technical but escential 
results of sections 3 and 4.

We would like to thank Prof. M. Shub for comments 
and suggestion about the stable manifold theorems, and 
Prof. G. Corach who told us  about the Aluthge transform, and 
shared with us fruitful discussions concerning these matters.

\section{Preliminaries.}
In this paper $\mat$ denotes the algebra of complex $r\times r$ matrices, 
$\matinv$ the group of all invertible elements of $\mat$, $\matu$ the group 
of unitary operators, and $\matsa$ (resp. $\matah$) denotes the real subspace 
of hermitian (resp. antihermitian) matrices. We denote 
$\cN(r) = \{ N \in \mat : N $ is normal$\}$. 
If $v \in \cene$, we denote by 
$\mbox{\rm diag}(v) \in \mat$ the diagonal matrix with $v$ in its diagonal. 

Given $T \in \mat$, $R(T)$ denotes the
range or image of $T$, $\ker(T)$ the null space of $T$, 
$\mbox{ rk}(T)= \dim R(T)$ the rank of $T$, $\sigma (T)$ the spectrum of $T$, 
$\la(T) \in \cene$ the vector of eigenvalues of $T$
  (counted with multiplicity), 
$\rho(T)$ the spectral radius of $T$, $\tr(T)$ the trace of $T$,
and $T^*$ the adjoint of $T$. 
We shall consider the space of matrices $\mat$ as a real Hilbert space 
with the inner product defined by
$$
\pint{A,\ B}=\Preal\big(\tr(B^*A)\big).
$$
The norm induced by this inner product is the 
Frobenius norm, that is denoted by $\|\cdot \|_2$. 
For $T \in \mat$ and $\cA\inc \mat$, by means of $\dist(T, \cA)$ we denote 
the distance between them, with respect to the Frobenius norm.

\medskip
On the other hand, let $M$ be a manifold. By means of $TM$ we denote the tangent bundle of $M$ 
and by means of $T_xM$ we denote the tangent space at the point $x\in M$. Given a function 
$f\in C^{k}(M)$ ($ k\ge 1$), 
we denote by $\der{f}{x}{V}$  the derivative of $f$ at the point $x$ applied to  
the tangent vector  $V \in T_xM$. 
Given an open subset $\cU\subseteq \R^m$, in $C^k(\cU,M)$ we shall consider the $C^k$-topology. 
In this topology,  two functions are close if the functions and all their derivatives until 
order $k$ are close uniformly on compact subsets. We denote by 
 $\mbox{Emb}^k(\cU ,M) $ the subset of $ C^k(\cU,M)$ consisting of the embeddings from 
 $\cU$ into $M$, endowed with the relative $C^k$-topology.

\subsection{Basic facts about the Aluthge transform}

\begin{fed}\rm
Let $T\in\mat$, and suppose that $T=U|T|$ is the polar
decomposition of $T$. Then,  the Aluthge transform of $T$ is defined by
\begin{align*}
\alu{T}&=\aluf{T}
\end{align*}
Throughout this paper, $\aluit{n}{T}$ denotes the n-times iterated Aluthge transform of $T$, i.e.
\begin{align*}
\aluit{0}{T}=T; \peso{and}
\aluit{n}{T}=\alu{\aluit{n-1}{T}}\quad n\in\N.
\end{align*}
\end{fed}

\noi The following proposition contains some properties of the Aluthge transform which follows easily from its definition.

\begin{pro}\label{facilongas} \rm
Let $T\in\mat$. Then:
\begin{enumerate}
  \item $\alu{T} = T $ \sii $T $ is normal. 
  \item $\alu{\la T}= \la \alu{T}$ for every $\la \in \C$.
  \item $\alu{VTV^*}=V\alu{T}V^*$ for every $V\in\matu$.
  \item $\|\alu{T}\|_2\leqp \|T\|_2\,$. In \cite{[AMS]} it is proved that equality only holds
  if $T \in \cN(r)$. 
  \item $T$ and $\alu{T}$ have the same  characteristic polynomial, in particular,
  $\spec{\alu{T}}=\spec{T}$.
\end{enumerate}
\end{pro}

\noi
The following result is easy to see, and it will be very useful in the sequel. 
\begin{pro}\label{sumas} \rm 
If $T=T_1\oplus T_2 \in \mat $ with respect to a reducing subspace 
$\ese \inc \cene$,  then   $\alu{T}=\alu{T_1}\oplus \alu{T_2}$.
  \end{pro}
\noi Next theorem states the regularity properties of Aluthge transforms 
(see \cite{[Dykema]} or \cite{[APS]}).

\begin{teo}\label{continuidad} \rm 
The map $\Delta$  is 
continuous in $\mat$ and it is of class $C^\infty$ in $\matinv$.
\end{teo}

\begin{rem} \label{fallo}\rm 
The map $\Delta$  fails to be differentiable at several matrices 
$T\in \mat \setminus \matinv $. Indeed, supose that $\Delta $ were differentiable
at $T = 0$. In this case, given $X \in \mat$, 
$$
T_0 \Delta (X) = \frac {d}{dt} \, \alu{tX}\Big|_{t=0} = \frac {d}{dt} \, t \, \alu{X}\Big|_{t=0} 
= \alu {X} \ .
$$
But this is impossible, because the map $X \mapsto \alu{X} $ is not linear. 
Using Proposition  \ref {sumas}, this fact
can be easily extended to any $T \in \mat\setminus \matinv$ such that $\ker T $ is 
orthogonal to $R(T)$ (for example, every non invertible normal matrix). 
\EOE
\end{rem}

\noi Now, we recall a result proved 
by Jung, Ko and Pearcy in \cite{[JKP1]}, and  by Ando in \cite{[Ando]}.

\begin{pro}\label{puntos limites normales} \rm 
If $T\in\mat$,  the limit points of the sequence
$\{\aluit{n}{T}\}_{n\in \N}$ are normal. Moreover, if $L$ is a limit point,
then $\spec{L}=\spec{T}$ with the same algebraic multiplicity. \QED
\end{pro}

\begin{rem}\label{uno solo} \rm 
Proposition \ref{puntos limites normales} has the next easy consequences: 
\ben
\item Let $T \in \mat$. Denote by $\la(T) \in \cene$ the vector of eigenvalues of $T$  
  (counted with multiplicity). Then 
 $\ds \|\aluit{n}{T}\|_2 ^2 \convd \sum _{i = 1}^r 
  |\la_i (T)|^2$.
\item If $\sigma (T) = \{\mu\}$, then $\aluit{n}{T} \conv \mu\,I$, because
$\mu \, I $ 
 is the unique normal matrix with spectrum $\{ \mu\}$. \EOE
 \een
\end{rem}

\subsection{Similarity orbits}
\begin{fed} \rm
Let  $D\in\mat$. By means of $\orb{D}$ we denote the similarity orbit of $D$:  
$$
\orb{D} = \{ \ SDS\inv\  : \ S \in \matinv \ \} \ . 
$$
On the other hand, $\orbu{D}= \{ \ UDU^* \  : \ U \in \matu \ \}$ 
denotes the unitary orbit of $D$. 
We donote by 
$\sub{\pi}{D} : \matinv \to \orb{D} \inc \mat $  the $C^\infty$ map defined by 
$\sub{\pi}{D}(S) = SDS\inv $, for every $S \in \matinv$.  
With the same name we note its restriction to the unitary group:  
$\sub{\pi}{D} : \matu \to \orbu{D} $. 
\EOE
\end{fed}
\begin{rem} \label{decrece}\rm
Let $T\in \mat$ and  $N \in \cN(r)$ with $\la(N) = \la(T)$. Then 
 $$
    \orbu{N}    =\{M \in \cN(r) : \la (M) = \la (T)\} \ .
    $$ 
  On the other hand, by Schur's triangulation theorem, there exists 
  $N_0 \in \orbu{N} $ such that 
  $\ds 
  \|T\|_{_2} ^2 - \sum _{i = 1}^r  |\la_i (T)|^2  = \|T-N_0\|_{_2}^2 \ge 
  \dist (T, \orbu{N}\, )^2  \ .
  $
  \EOE
\end{rem}

\noi 
The following two results are well known (see, for example, \cite{[CPR]} or \cite{[AS]}):
\begin{pro}\label{son variedades}\rm
The similarity orbit $\orb{D}$ is a $C^\infty$ submanifold of $\mat$, and the projection
$\sub{\pi}{D} : \matinv \to \orb{D}$ becomes a submersion. Moreover, $\orbu{D}$ is a compact submanifold of $\orb{D}$, which consists of the normal elements of $\orb{D}$, and 
$\sub{\pi}{D} : \matu \to \orbu{D}$ is a submersion.\QED
\end{pro}

\begin{rem} \label{tan} \rm
For every $N\in\orb{D}$, it is well known that
\begin{align}
\sub{T}{N}\,\orb{D}&= \sub{T}{I}( \pi_N ) (\mat\,)  =\{[A,N]=AN-NA: \ A\in\mat\}. \nonumber\\
\intertext{If $\spec{D } = \{ \mu_1 , \dots , \mu_k \} $, and $E_i(N) $ 
are the spectral projections of $N \in \orb{D}$ associated to disjoint open neighborhoods of 
each $\mu_i\,$, then $N = \sum_{i=1 }^ k \mu_i E_i (N)$. Therefore  } 
\sub{T}{N}\,\orb{D}& =\{AN-NA: \ A\in\mat\} \nonumber\\ & = \left\{ \sum _{i, j = 1}^k (\mu_j - \mu_i ) 
E_i(N) A E_j(N): \ A\in\mat \right\}\nonumber \\
\label{TUN} & = \{ X \in \mat : E_i(N) X E_i(N) = 0 \,,\ 1\le i \le k\} \ . 
\end{align}
Throughout this paper we shall consider on $\orb{D}$ (and in $\orbu{D}\,$) the 
Riemannian structure inherited from $\mat$ (using the usual inner product on their 
tangent spaces). 

Observe that, for every $U \in \matu$, it holds that $U\orb{D} U^* = \orb{D} $ and the map 
$T \mapsto UTU^*$ is isometric, on $\orb{D}$,  with  respect to the Riemannian metric as 
well as with respect to the $\| \cdot \|_2$ metric of $\mat$. \EOE
\end{rem}

\medskip

\noi
In the previous work \cite {[APS]}, we have proved the following result:

\begin{teo}\label{the key} \rm 
Let $D= \mbox{\rm diag}(d_1,\ldots,d_r)  \in \mat$ be an invertible diagonal matrix.
For every $N\in\orbu{D}$, there exists a 
subspace $\sub{\cE}{N}^s$ of the tangent space $\sub{T}{N}\orb{D}$ such that
\begin{enumerate}
	\item $\sub{T}{N}\orb{D}=\sub{\cE}{N}^s\oplus \sub{T}{N}\orbu{D}$;
	\item Both, $\sub{\cE}{N}^s$ and  $\sub{T}{N}\orbu{D}$, are $T\,\Delta$-invariant;
	\item $T\,\Delta_N\big|_{\sub{T}{N}\orbu{D}} = I_{\sub{T}{N}\orbu{D}}$ and  
	$\left\|T\,\Delta_N\big|_{\sub{\cE}{N}^s}\right\|\leq\sub{k}{D}<1$, where 
	$$
	\sub{k}{D}=\max_{i,\,j\,:\ d_i\neq d_j}\frac{|1+e^{i(\arg(d_j)-\arg(d_i))}|\,  |d_i|^{1/2}|d_j|^{1/2}}{|d_i|+|d_j|}\,;
	$$
  \item If $U\in\matu$ satisfies $N=UDU^*$, then $\sub{\cE}{N}^s=U(\sub{\cE}{D}^s)U^*$.
  \end{enumerate}
  In particular, the map $\orbu{D} \ni N \mapsto \sub{\cE}{N}^s$ is smooth. This fact can be 
  formulated in terms of the projections $P_N$ onto $\sub{\cE}{N}^s$ parallel to $\sub{T}{N}\orbu{D}$, 
  $N \in \orbu{D}$. \QED
\end{teo}

\newcommand{\diag}[1]{\hbox{\rm diag}\left( #1\right)}

\subsection{Stable manifold theorem for the Basin of attraction}

Let $M$ be a smooth Riemann manifold, 
$f$ a smooth endomorphism of $M$, and 
$N\subseteq M$ a compact set  such that $f(N)=N$. 
The basin of attraction of $N$ is the set 
$$
B_f(N)=\{y\in M: \dist(f^n(y),N)\conv 0 \} \ .
$$
Given $\eps>0$,  the local basin of $N$ is the set
$$
B_f(N)_\eps = \{y\in B_f(N) :  \dist(f^n(y),N)<\eps \ ,\ \  \mbox{for every $n \in \N$} \} \ .
$$
The following result 
is standard when stated on a compact $f$-invariant subset $N$ 
 (see the Appendix of \cite{[APS]}, \cite {[HPS]} or \cite{Sh2}).  
The following version, which extend the prelamination $\ewe^s$ to its local basin  $B_f(N)_\eps$
is also well known. In Remark \ref{esquech}, we shall expose briefly the principal steps of the
proof of the version for $N$, and then explain how that proof can be
extended to ``its basin of attraction."

\begin{teo}[Stable manifold theorem]\label{teorema 5.5} \rm
Let $f$ be a $C^k$ endomorphism of $M$ and let $N$ be a compact $f$-invariant subset of $M$.  
Let us assume that for some $\eps>0$ there exist  
two continuous subbundles of $\sub{T}{B_f(N)_\eps}M$, denoted by $\cE^s$ and $\cF$, 
such that, for every $x \in  B_f(N)_\eps \,$, 
\begin{enumerate}
	\item $\sub{T}{B_f(N)_\eps}M= \cE^s \oplus \cF$.
	\item $\cE_x^s$ 
	is $T_xf$-invariant in the sense that $T_x f( \cE_x^s ) \inc  \cE_{f(x)}^s \,$. 
	\item $\cF_z$ is $T_z f$-invariant, for every $z \in N$. 
	\item There exists $\rho\in (0,1)$ such that $T_x\, f$ restricted to $\cF_x$  
	expand it by a factor 	greater than $\rho$,  
	and $\dersin{f}{x}:\cE_x^s\to\cE_{f(x)}^s$ has norm lower than $\rho$. 
	
\end{enumerate}
Then, there is a continuous, $f$-invariant and self coherent 
$C^k$-pre-lamination  $\ewe^s: B_f(N)_\eps \to \mbox{Emb}^k((-1,1)^m,M)$ 
(endowed with the $C^k$-topology) such that, for every $x \in  B_f(N)_\eps \,$, 
\begin{enumerate}
\item $\ewe^s(x)(0)=x$,
\item  $\ewe_x^s=\ewe^s(x)((-1,1)^m)$  is tangent to $\cE_x^s \,$, 
\item $
\ewe_x^s \inc \Big\{y\in M\ : \  \dist(f^n(x),f^n(y))<\dist(x,y)\rho^n\Big\}
$. \QED
\end{enumerate}
\end{teo}

\begin{rem} \label{esquech}\rm 
The proof of the existence of a map $\ewe^s: N \to \mbox{Emb}^k((-1,1)^m,M)$ 
which satisfies all the the mentioned conditions, consists in using the graph transform 
operator.
We shall see that it is well defined if we only consider forward iterates. 
Therefore, since the basin of attraction of $N$ is properly mapped
inside by $f$,  the graph transform operator is well defined on
$B_f(N)_\eps\,$, allowing to extend the proof of stable manifolds to the whole
local basin.
Recall that to define the graph transform operator, first we consider  
 $
 C^k(\hat{\cE^s}_x,\hat{\cF}_x)
 $,  
 the set of $C^k$ maps from $\hat{\cE^s}_x$ to $\hat{\cF}_x\,$, where 
 $$
 \hat{\cE^s}_x(\mu)= \exp(\cE^s_x\cap {(T_xM)}_\mu),\,\,\,\,\, \hat{\cF}_x(\mu)= \exp(\cF_x\cap {(T_xM)}_\mu)
 $$ 
 and $\exp_x: ({T_xM})_\mu\to M$ is  the exponential map acting on $({T_xM})_\mu\,$, 
 the ball of radius $\mu$ in $T_xM$. Later we  consider the space 
$$C^{k,0}(\hat{\cE^s},\hat{\cF})= \{\sigma: N\to C^k(\hat{\cE^s}_x,\hat{\cF}_x)\}$$ i.e.: 
for each $x\in N$ we take  $\si_x\in C^k(\hat{\cE^s}_x,\hat{\cF}_x)$ and we assume that 
$x\mapsto \si_x$ moves continuously with $x$. We can represent $C^{k,0}(\hat{\cE^s},\hat{\cF})$ 
as a vector bundle over $N$ given by $N\times \{C^k(\hat{\cE^s}_x, \hat{\cF})\}_{x\in X}\,$. 
Then, we take  the maps 
$$ 
f^1_x= p^1_x\circ f: M\to \hat{\cE^s}_x \peso{and}  f^2_x= p^2_x\circ f: M\to \hat{\cF}_x \ ,
$$ 
where $p_x^1$ is the projection on $\hat{\cE^s}_x$ and $p_x^2$ is the projection on $\hat{\cF}_x\,$.
Now we take the {\it graph transform operator}. If $f$ is a diffeomorphism, then we can 
obtain an explicit formula for the graph transform: 
\begin{eqnarray}\label{gto}
\Gamma_f(\sigma_x)= \big(f_x^2\circ(id,\si_{f(x)})\,\big)^{-1}
\circ \big(f_x^1\circ(id,\si_x)\,\big)\big | _{\hat{\cE^s}_x} \ .
\end{eqnarray}
On the other hand, if $f$ is an endomorphism, the graph transform can be defined implicitly.
In both cases, this map is well defined in $B_f(N)_\eps$ and therefore the whole proof 
can be carried out, in the sense of  proving that the graph transform operator is a 
contractive map and therefore it has a fixed point. \EOE
\end{rem}

\subsection {Spectral projections}
In this section we state the basic properties 
of the spectral projections of matrices, which are constructed by using the
Riesz functional calculus. A complete exposition on this theory can be found in 
Kato's book \cite[Ch. 2]{Kato}.
Given  $T \in \mat$ we call $\la = \la (T) \in \cene $ its vector of eigenvalues. 
Let $\spec{T } = \{ \mu_1 , \dots , \mu_k \} $, taking one $\mu_i$ for 
each group of repeated  $\la_j(T) = \mu_i$ in $\la (T)$ (i.e., $k \le r$). Fix $D= \diag{\la} \in \mat $.  

\begin{fed}\label{laymu}\rm
Given $\la = \la (D)\in \cene$ and  
$\mu = ( \mu_1 , \dots , \mu_k) \in \C^k$ as before ($\mu_i \neq \mu_j$), let 
\ben
\item $\ds \eps_\mu = \frac13 \, \min_{i \neq j} |\mu_i - \mu_j|\ $ 
and $\ds \Omega_\mu = \bigcup _{1\le i \le k} B(\mu_i \, , \, \eps_\mu)$. 
\item $\tilde { \cM}_{\mu } = 
\{ T \in \mat : \spec{T} \inc \Omega_\mu 
 \}$, which is open in  $\mat$.
\item Let $E : \tilde { \cM}_\mu \to \mat ^k  $ given by 
$$
\tilde { \cM}_\mu \ni T \mapsto E(T) = (E_1(T) , \cdots , E_k(T)\,)\,,
$$ 
where $E_i (T) = \aleph_{B(\mu_i \, , \, \eps_\mu)}(T) $ is the 
spectral projection of  
$T\in \tilde { \cM}_\mu\,$,  associated to $B(\mu_i \, , \, \eps_\mu)$ . 
\item Denote $Q = (Q_1, \dots , Q_k ) = E(D)$ and consider the open set 
 \beq\label{ml}
 \cM_\la = \{ T \in \tilde { \cM} _\mu : \mbox{ rk}(E_i (T) ) = \mbox{ rk}(Q_i) \ ,  \ \ 1\le i \le k\} \ ,
 \end{equation}
which is  the  connected component of $D$  in $\tilde { \cM}_{\mu } \,$. 
 \item 
Let $\Pi _ E : \cM_\la \to \mat$ given by 
$\ds \Pi_E (T) = \sum_{i=1}^k \mu_i E_i (T)$, for every $T \in \cM_\la \,$.  
\EOE
 \een
\end{fed}

\begin{rem} \label{E y PIE} \rm  
Given $\la = \la (D)\in \cene$ and 
$\mu = ( \mu_1 , \dots , \mu_k) \in \C^k$ as before,  the  following properties hold:
\ben 
\item  For every $1 \le i \le k$, $Q_i = Q_i ^*\,$. Also  
 $Q_i Q_j = 0$ (if $i\neq j$) and $\ds \sum_i Q_i = I$. 
The entries of  $E(T)$ for other  $T \in \cM_\la$ satisfy the same properties, but they can be not selfadjoint. 
\item Each  map $E_i$ (so that the map $E$)  is of class $C^\infty$ in $\cM_\la \,$.
\item $E(\cM_\la ) = \orb {Q}: = \{ 
(SQ_1S\inv , \dots , SQ_kS\inv ) : S \in \matinv \} $. 
\item Moreover, if  $T \in \cM_\la $ y $S \in \matinv$, 
then  $E(STS\inv ) = S E(T) S\inv$. 
\een
Then,  the map $\Pi _ E : \cM_\la \to \mat$ 
satisfies the following properties: 
\ben 
\item It 
is of class  $C^\infty$ on $\cM_\la \,$.
\item For every $T \in \orb{D}$, $\ds \Pi_E (T) = T$. 
\item $\ds \Pi_E (\cM_\la ) = \orb {D} $, and $\rho(T - \Pi_E (T)\,) < \eps_\mu$ for every $T \in \cM_\la \,$.  
\item If $T \in \cM_\la $ and $S \in \matinv$, then 
 $\ds \Pi_E (STS\inv ) = S  \Pi_E (T) S\inv$. \EOE
\een
\end{rem}

\begin{rem} \label{AM} \rm
With the previous notations, for every $M\in \cM_\la\,$, we consider the subspace 
\beq\label{am}
  \cA_M  = \{ B \in \mat : BE_i (M) = E_i(M) B \ \ \ 1\le i \le k\}  \ ,
 \end{equation}
 of $block$ $diagonal$ matrices, with respect to $E(M)$. It is easy to see that 
$\cA_M = \ker T_M \Pi_E \,$ and $R (T_M \Pi_E ) =  T_{N}\orb{D}$.  By Eq. \eqref{TUN},  if $N = \Pi_E (M)\in \orb{D} $, then 
$\mat = \cA_M \oplus T_{N}\orb{D}$, and the sum becomes orthogonal if $M \in \orbu{D}$. 

Since $\Pi_E\, ^2 = \Pi_E\,$, if $M \in \orb{D}$, then $T_M \Pi_E  $ is the 
projector with kernel $\cA_M$ and  and image
$
T_{M}\orb{D}\,$. Observe that
\beq\label{am2}
\cA_M = \{M\}' :=  \{B \in \mat : MB = BM\}
 \end{equation}
for every $M \in \orb{D}$, since in this case $M = \Pi_E(M)$. \EOE

\end{rem}

\section{Some dynamical aspects of the Aluthge transform}

The main aim of this section is to introduce the dynamical setting as well 
as some results that will be used in the next section to prove the 
convergence of the iterated Aluthge transform sequence.

\subsection {The derivative of $\Delta$ in $\mat$}
Let $N\in\mat$ be an invertible normal matrix. Theorem \ref{the key} gives 
a description of the action of $T_N\Delta$  on 
$T_N \orb{N}$. By Remark \ref{tan}, in order to obtain a complete characterization 
of the action of $T_N\Delta$ on $\mat$, it is enough to describe the action of 
$T_N\Delta$ on its orthogonal complement, i.e. the subspace $\cA_N$ described in Remark 
\ref{AM}. 

\begin{pro} \label{foto}
Let $D= \mbox{\rm diag}(d_1,\ldots,d_r)  \in \matinv$ be a diagonal matrix
with $k$ different eigenvectors. Fix $N\in \orbu{D}$ and consider the subspace 
$\cA_N\inc \mat $ defined in Eq. \eqref{am}. Then $T_{N}\Delta  |_{\cA_{N}} = I_{\cA_{N}} \,$. 
\end{pro}
\bdem If $N \in \orbu{D}$, then $N$ is normal, and $\cA_{N} = \{N\}'$. 
Hence, if $Y \in \cA_{N}$ is normal, then $N+t Y $ is also normal for every 
$t \in \R$. This implies that 
$T_N\Delta  (Y) = \frac{d}{dt}\alu{N + tY} \Big|_{t=0}  = Y $.
On the other hand, since $\cA_N$ is closed by taking adjoints, 
then $\R e (X) \in \cA_{N}$ and $\mathbb I m (X) \in \cA_{N}$
for every $X \in \cA_{N}\,$. Therefore 
$T_N\Delta  (X) = X $.
\edem

\begin{cor} \label{mitades} 
Let $D= \mbox{\rm diag}(d_1,\ldots,d_r)  \in \mat$ be an invertible diagonal matrix. 
For $N\in \orbu{D}$, denote by $F_N = \cA _N \oplus T_N\orbu{D}\,$. Then, the subspaces 
$F_N$ and $\sub{\cE}{N}^s$ (defined in Theorem \ref{the key}) satisfy
that,  for every $N \in \orbu{D}$, 
\ben
\item $\mat = F_N \oplus \sub{\cE}{N}^s\,$. 
\item Both subspaces are $T_N\Delta $ invariant. 
\item $T_N\Delta \big|_{F_N} = I_{F_N}$ and 
$\left\|\left. T_N\,\Delta\right|_{\sub{\cE}{N}^s}\right\|\leq\sub{k}{D}<1$, 
where $ \sub{k}{D}$ is defined as in Theorem \ref{the key}.
\item The distributions $N \mapsto F_N$ and $N \mapsto \sub{\cE}{N}^s$  are smooth. 
\een
\end{cor}
\bdem By Theorem \ref{the key} and  Remark \ref{AM}, 
$$
\mat =  \cA_N \oplus T_N \orb{D}  = \cA_N \oplus T_N \orbu{D} \oplus \sub{\cE}{N}^s 
=  F_N \oplus \sub{\cE}{N}^s \ . 
$$
By Proposition \ref{foto}, one deduces that  $T_N\Delta \big|_{F_N} = I_{F_N}\,$. 
The remainder conditions follow easily from  Theorem \ref{the key}.
\edem
 
 \begin{rem} \label{CFD} \rm
 With the notations of Corollary \ref{mitades}, the 
 subspaces $F_N$ and $\sub{\cE}{N}^s$ can be characterized by mean of the 
 functional calculus applied to the linear maps $T_{N} \Delta $, for every 
 $N \in \orbu{D}$.  Indeed,  
 $$ F_N = R\big( \aleph _{B(1 , \eps )} (T_{N} \Delta  \,) \,\big)
 \peso{and}  
 \sub{\cE}{N}^s = R\big( \aleph _{B(0 , k_D + \eps )} (T_{N} \Delta  \,) \, )\ , 
 $$
for every  $\eps >0$ sufficiently small. In particular, this implies that the 
distribution of subspaces  $F_N$ and $\sub{\cE}{N}^s$ can be extended smoothly 
to an open neighborhood of $\orbu{D}$.\EOE
\end{rem}

\subsection{Stable manifolds}
Let $\cP \inc \mat$ be a compact set of fixed points for $\Delta$, i.e. a compact set of normal matrices. 
Recall that 
its basin of attraction is the set 
$$
B_\Delta (\cP) = \{ T \in \mat : \dist (\aluit{n}{T}, \cP)  \conv 0 \} \ 
$$
and, for every   $\eps > 0$, the local basin  is the set 
$$B_\Delta (\cP)_\eps = 
\{T \in B_\Delta (\cP) : \dist (\aluit{n}{T} , \cP) <\eps\, , \ n \in \N\} \ .
$$ 
In this subsection, using the stable manifold theorem \ref{teorema 5.5}, we shall prove that, 
if $\cP$ has a distribution of subspaces with good properties 
(like the distribution of Corollary \ref{mitades} for $\cP = \orbu{D}\,$),  
through each $T\in B_\Delta(\cP)$ closed enough to $\cP$ there 
is a stable manifold $\ewe^{ss}_T$ with the property
$$
\ewe^{ss}_{T} \subseteq \{B\in\mat: \ \|\aluit{n}{T} - \aluit{n}{B} \|< C\gamma ^n \},
$$ 
where $\gamma<1$ and $C$ is a positive constant. With this aim, firstly we need 
to extend to some local basin the distribution of subspaces given on $\cP$.  
This extension is a quite standard procedure in dynamical systems. 
For completeness, we include a sketch of its proof 
(adapted to our case) in the Appendix A. 

\begin{pro}\label{conitosA} 
Let $\cP $ be a compact set consisting of fixed points of $\Delta$. 
Suppose that, for every $N \in \cP$, 
there are subspaces 
$\sub{\cE}{N}^s$ and $\sub{\cF}{N}$ of $\mat$ 
with the following properties:
\begin{enumerate}
	\item $\mat=\sub{\cE}{N}^s\oplus \sub{\cF}{N}\,$.
			\item The distributions $N \mapsto \sub{\cE}{N}^s $ and $N \mapsto \sub{\cF}{N}$ are continuous. 
		\item There exist $\rho  \in(0,1)$ which does not depend on $N$ such that 
	\begin{equation}\label{menor que rho}
	\left\|\left. T_N\,\Delta\right|_{\sub{\cE}{N}^s}\right\| < 1- \rho \peso{and} 
\left\|\left. (I - 	T_N \,\Delta\right)|_{\sub{\cF}{N}}\right\|<\frac \rho2 \ ,
	\end{equation}
\item Both subspaces $\sub{\cE}{N}^s$ and $\sub{\cF}{N}$ are $T_N\Delta $ invariant. 
\een
Then, there exists $\eps > 0$  such that 
the distributions $N \mapsto \sub{\cE}{N}^s $ and $N \mapsto \sub{\cF}{N}$ can be extended 
to the local basin 
$B_\Delta (\cP)_\eps\,$,  
verifying conditions 1, 2, 3 and the following new condition: 

\ben
\item [4'] For every $T \in B_\Delta (\cP)_\eps\,$, the subspace $\cE_T^s$ is  $T_T \, \Delta$-invariant, i.e., 
 $
  \sub{T}{T}\Delta(\sub{\cE}{T}^s)\subseteq \sub{\cE}{\alu{T}}^s  \ .
  $
\end{enumerate}
\end{pro}

\bdem See section A.1 of the appendix. 

\bigskip

\noi  Now we are ready to state and prove the announced result on stable manifolds.

\begin{pro}
\label{pelitosA}
Let $\cP $ be a compact set consisting of fixed points of $\Delta$, with two distributions 
$N \mapsto \sub{\cE}{N}^s $ and $N \mapsto \sub{\cF}{N}$ 
which satisfy the hypothesis of 
Proposition \ref{conitosA}. 
Then, 
there exist $\eps>0$  and a  $C^2$-pre-lamination \rm $\ewe^s: B_\Delta (\cP)_\eps  \ \to 
\mbox{Emb}^2((-1,1)^m,B_\Delta (\cP) )$ \it (endowed with the $C^2$-topology) 
of class $C^0$ such that, 
 for every $T\in B_\Delta (\cP)_\eps \,$, 
\begin{enumerate}
  \item  $\ewe^s (T) (0) = T$. 
  \item If $\ewe^{ss}_T $ is 
 the submanifold $\ewe^s (T) \big( (-1,1)^m \big)$, then $\sub{T}{T}\ewe^{ss}_T=\sub{\cE}{T}^s\,$.
	\item There are constants $\gamma<1$ and $C>0$ such that 
\beq\label{lasW}
\ewe^{ss}_{T} \subseteq \{B\in\mat: \ \|\aluit{n}{T} - \aluit{n}{B} \|< C\gamma ^n \} .
\end{equation}
\end{enumerate}
\end{pro}

\bdem By Proposition \ref{conitosA}, the distributions $\cP \ni N \mapsto \sub{\cE}{N}^s \, , \  \sub{\cF}{N}$
can be extended to a local basin $B_\Delta (\cP)_\eps \,$, satisfying  the hypothesis 
the Stable Manifold Theorem \ref{teorema 5.5}. 
Observe that the condition 
$\left\|\left. (I - 	T_T \,\Delta\right)|_{\sub{\cF}{T}}\right\|<\frac \rho2 $
implies that $\|T_T \,\Delta (Y)\| >  (1- \frac\rho2) \|Y\| $ for every $
Y \in \sub{\cF}{T} \,$.  
\edem

\section{Convergence of the 
sequence $\aluit{n}{T} $}
This section is entirely devoted to the proof of Jung, Ko and Pearcy's conjecture. The basic tools 
are the results of the previous section.

\subsection{The case $\cP = \orbu{D}$}
Let $D\in\mat$ is an invertible diagonal matrix. In this section we shall 
consider the compact invariant set $\cP=\orbu{D}$.    
Observe that the distributions $N \mapsto F_N$ and $N \mapsto \sub{\cE}{N}^s$  
($N \in \cP$) given  by Corollary \ref{mitades} clearly verify the hypothesis of 
Propositions \ref{conitosA} and \ref{pelitosA}. 
Thus, by Proposition \ref{pelitosA}, there exist a continuous prelamination 
$\ewe^s: B_\Delta (\cP)_\eps \to \mbox{\rm Emb}^2((-1,1)^k,B_\Delta (\cP))$  
and submanifolds $\ewe^{ss}_T $ (for every $T \in B_\Delta (\cP)_\eps$). 
which  will be also used throughout this section. 

In this setting we can give simple characterizations  
of the basins of $\cP$. Indeed,  
by Proposition \ref{puntos limites normales} and Remark \ref{decrece}, 
\beq\label{bd1}
B_\Delta (\orbu{D}\,)= \{ T \in \mat :  \la (T)  = \la (D) \}\ .
\end{equation}
Given $T \in  B_\Delta (\orbu{D}\,)$, let
$d_n(T) =    \|\aluit{n}{T} \|_{_2} ^2 - \sum _{i = 1}^r  |\la_i (T)|^2$. 
By Remarks \ref{uno solo} and \ref{decrece}, 
$\ds \dist (\aluit{n}{T} ,  \orbu{N} \,)  \le d_n(T)   \convd 0 \ .
$
Then, 
\beq\label{bd2}
\{ T \in  B_\Delta (\orbu{D}\,) : d_1(T) < \eps \}
\inc B_\Delta (\orbu{D}\,)_\eps \ , 
\end{equation}
and it is also an open neighborhood of $\orbu{D}$ in $ B_\Delta (\orbu{D}\,)$. Therefore, 
if $T\in B_\Delta (\orbu{D}$ is close enough to $\orbu{D}$, then 
$T\in B_\Delta (\orbu{D}\,)_\eps$ (and we do not need to check the $\dist (\aluit{n}{T} ,  \orbu{N} \,) $ 
for $n >1$).

\medskip

\noi Observe that if $T \in B_\Delta (\orbu{D}\,)$, despite the equality $\la (T) = \la (D)$, $T$ can have any 
Jordan form.

\subsection{The sets $\lentas$}

In this subsection we identify some convenient sets of matrices where 
the iterated Aluthge transform sequence 
converges (possibly slowly). 
By their properties, these sets will play a key role in the proof of the convergence of 
the iterated Aluthge transform sequence. Let $D$ be an invertible diagonal matrix, 
 $\la = \la(D)$, and  
$\Pi_E : \cM_\la \to \orb{D}$, the map defined in Section 2.4.  
If $\cP = \orbu{D}$, consider the following subset of $B_\Delta(\cP)$:
\beq\label{o de d}
\lentas = \{T \in B_\Delta(\cP) : \Pi_E (T) \in \cP \}=\Pi_E^{-1}(\cP) \cap B_\Delta(\cP) \ .
\end{equation} 
Note that, if  $T \in \cO_D\,$,  then the system of projectors $E(T)$ is orthogonal. 
Hence we get the next simple consequence.

\begin{pro}\label{en oded converge} 
If $T \in \lentas \,$, then $\aluit{n}{T} \conv  \Pi_E (T) \in \orbu{D}$. 
\end{pro}
\bdem
If $T \in \lentas\inc  B_\Delta(\orbu{D}\,)$ then $\la (T) = \la (D)$, 
by Eq.  \eqref {bd1}. On the other hand, if  $N = \Pi_E (T)$, 
then $E(T) = E(N)$ is an {\bf orthogonal} system of projectors, and
$T \in \cA_N\,$, the subspace defined in Eq. \eqref{am}.  
Write $T = T_1 \oplus \dots \oplus T_k\,$, where each $T_i = T|_{R(E_i(T)\,)}\,$. 
By Proposition \ref{sumas},  $\Delta ^n (T) = 
 \Delta ^n (T_1 ) \oplus \dots \oplus \Delta ^n (T_k )\,$, for every $n \in \N$. Since 
 $\la (T) = \la (D)$,  then  $\sigma (T_i ) = \{\mu_i\}$ and, 
 by Remark \ref{uno solo}, 
 $$ 
 \aluit{n}{T_i } \conv \mu_i \, I_{R(E_i(T)\,)} \ ,
 \ \ 1\le i \le k \ . 
 $$
 
Therefore $\ds \Delta ^n (T) \conv \sum_{i=1}^k \mu_i \, E_i(T) = \Pi_E (T)  \ .
$ 
\edem

\noi Another important characteristic of the sets $\lentas$ is that each 
element of $B_\Delta (\cP)$ ``close enough" to $\cP$ is exponentially attracted toward $\lentas \,$. This property is precisely described in the following statement. 

\begin{pro}\label{los w cortan} \rm 
Let $D\in\mat$ be an invertible diagonal matrix, $\cP=\orbu{D}$ and  
$\ewe^s: B_\Delta (\cP)_\eps \to \mbox{\rm Emb}^2((-1,1)^k,B_\Delta (\cP) )$  
the prelamination given by Proposition \ref{pelitosA}. 
Then, there exists $\eta<\eps$ such that $\ewe^{ss}_{ T} \cap \lentas \neq \varnothing$ for every 
$T \in  B_\Delta (\cP)_\eta \, $. 
\end{pro}

\noi
The proof is rather technical. 
Since this result is the key part of the proof of the conjecture, 
we give a brief description of the proof here, and we leave a detailed proof to Appendix A. 

\medskip
If $T\in B_\Delta (\cP)_\eps$ is near $N \in \orbu{D}$, then the set
$\cV^{ss}_{ T} = \Pi_E (\ewe^{ss}_{ T} \,)$ is a smooth submanifold of $\orb{D}$
with  the same dimension as $\ewe^{ss}_{ T} \,$, 
and it remains $C^2$-close to $\ewe^{ss}_{ N}\,$. This facts can be deduced from the properties 
of the projection $\Pi_E $ stated in Subsection 2.4. 

In order to show that 
$\ewe^{ss}_{ T} \cap \lentas \neq \varnothing$, it suffices to prove that 
$\cV^{ss}_{ T} \cap \orbu{D} \neq \varnothing$. Finally, this fact follows from a well known
argument of transversal intersection (inside the manifold $\orb{D}\,$), by using 
the dimension of the tangent spaces, and the fact that 
$\ewe^{ss}_{ N} $ intersects $\orbu{D}$ transversally.

\subsection{The proof of Jung, Ko and Pearcy's conjecture}

Now, we are in conditions 
to prove the main result of this paper:
\begin{teo} \label{EL TEO}
For every $T \in \mat$, the sequence $\aluit{n}{T} $ converges. 
\end{teo}
\bdem 
By Corollary 4.16 of \cite{[AMS]}, we can assume that $T \in \matinv$. Let $D\in  \matinv$ be a diagonal matrix 
such that $\la (T) = \la (D)$, and let $\cP = \orbu{D}$. 
By Eq. \eqref{bd1}, $T \in B_\Delta (\cP)$. By Eq. \eqref{bd2},  
replacing $T$ by  $\aluit{n}{T} $ for some $n$ large enough, we can assume that 
$T \in B_\Delta (\cP)_\rho \, $, for any fixed $\rho >0$. 

Consider now  the 
stable manifold $\ewe^{ss}_{ T}$,  constructed in 
Proposition \ref{pelitosA}, for every $T \in B_\Delta (\cP)_\eps \, $. 
By Proposition \ref{los w cortan}, there exists $0<\eta <\eps$ such that 
$\ewe^{ss}_{ T} \cap \lentas \neq \varnothing$, 
for every $T \in  B_\Delta (\cP)_\eta \, $. If $M \in \ewe^{ss}_{ T} \cap \lentas \,$ then, 
by  Proposition \ref {en oded converge} and Eq. \eqref{lasW} of Proposition \ref {pelitosA},  we deduce that 
$\ds \Pi_E (M) = \lim _{n\to \infty} \aluit{n}{M} 
= \lim _{n\to \infty} \aluit{n}{T} $.  
\edem

\section{Regularity of the map $\Delta^\infty$}
Given $T \in \mat$ we denote $\ds \aluit{\infty}{T} = \lim_{n \to \infty}  \aluit{n}{T} $, which is 
a normal matrix. 
Observe that the map $ \ds \Delta^{\infty}: \mat \to \cN(r)$ is a retraction. In this section 
we study the regularity of this retraction. 

\subsection{Differentiability vs. continuity}

\noi
In \cite{[APS]} we proved that the map $\Delta^\infty$ is of class $C^\infty$, when it is
restricted to the open  dense set of those matrices in $\mat$ with $r$ different eigenvalues. 
The following proposition shows that this can not be extended globally to the set of all matrices.

\begin{pro}\label{no C1}
The map $\Delta^\infty$ can not be $C^1$ in a neighborhood of  the identity.
\end{pro}
\bdem
Suppose that $\Delta^\infty$ is $C^1$ in a neighborhood of  the identity. 
By the same argument used in the proof of Proposition \ref{foto}, it follows that 
$\sub{T}{I}{\Delta^\infty}$ is the identity map 
 (in this case, $\cA_I = \mat$\,). 
This implies that $\Delta^\infty$ is a local diffeomorphism.  
However, this is impossible because it takes values in the set of normal operators.
\edem

\subsection
{Continuity of $\Delta^\infty$ on $\matinv$}
For a sake of convenience, throughout this subsection we shall use the spectral norm, instead of 
the Frobenius norm,  to measure distances in $\mat$.

\begin{rem}\label{reduccion} \rm 
Since $\Delta^\infty$ is the limit of continuous maps and it is a retraction, 
in order to show that  it  is continuous 
on $\matinv$, it is enough to prove the continuity at the normal matrices of $\matinv$. 
Indeed, observe that  $\Delta ^n$ is continuous for every $n \in \N$. 
Then, for every $T \in \matinv$, and every neighborhood $\cW$ of $\aluit{\infty}{T}$, 
there exists $n \in \N$ and a neighborhood $\cU$ of $T$ such that 
$\aluit{n}{\cU}\inc \cW$. Then, note that $\Delta^\infty \circ \Delta^n = \Delta^\infty$. 
\EOE  
\end{rem}

\noi
From now on, let $N_0 \in \cN(r)$ be a fixed normal invertible matrix such that 
$\la (N_0) = \la$ and $\spec{N_0}=(\mu_1,\ldots,\mu_k)$. 
Let $\eps_\mu = \frac13 \min_{i \neq j} |\mu_i - \mu_j|$. Consider the open set 
 $\cM_\la$ defined in equation \eqref{ml}. 
Recall that, for $T \in \cM_\la\,$, we call  $\ds \Pi_E (T) = \sum_{i=1}^k \mu_i E_i (T) 
\in \orb{N_0}\,$.  

\begin{fed}\label{pebeta} \rm With the previous notations, we denote  
\ben
\item 
$\mla{\eta}$  the open subset of $\cM_\la$ obtained in the same way, but by replacing $\eps_\mu$ by  
$\eta\in (0,\eps_\mu)$. Observe that $\rho(T - \Pi_E (T)\,) < \eta$ for every $T \in \mla{\eta} \,$.  
\item Given $\beta >0$,  we denote $\cP_\beta = 
\{ N \in \cN(r) : \dist (N , \orbu{N_0}\,) \le \beta\}\,$. Observe that 
$\cP_\beta $ is compact. 
\EOE
\een
\end{fed}

\begin{lem} \label{cosas} 
With the previous notations, let $\beta >0$ such that the closed ball $\overline{B(N_0 \, , \beta )}$ is contained in $ \cM_\la \,$. 
Then 
\ben
\item $\cP_\beta \inc \cM_\la \, $. 
\item For every $\eta < \min \{\beta , \eps_\mu \}$, it holds that $\mla{\eta} \inc B_\Delta (\cP_\beta )$. 
\item Moreover, if $\ T \in \mla{\eta}$ and $N = \aluit{\infty}{T}$, then $N \in \mla{\eta}\,$ and 
$\|N -\Pi_E (N) \| < \eta$.
\een
\end{lem} 
\bdem 
If $N \in \cP_\beta $, let $U \in \matu$ such that $\|N - UN_0 U^* \| \le \beta$ 
(recall that $\orbu{N_0}$ is compact). Then $U^* N U \in 
\overline{B(N_0 \, , \beta )} \inc \cM_\la \,$, so that also $N\in \cM_\la\,$.

Let $T \in \mla{\eta}$ and denote $N = \aluit{\infty}{T}$. Since 
$\la (N) = \la (T)$, then also $N \in \mla{\eta}\,$. 
  Since $N$ is normal, then $\Pi_E (N) \in \orbu{N_0}$ and 
  $$\eta > \rho ( N -\Pi_E (N) \,) = \|N -\Pi_E (N) \| \ ,$$ 
  because $N$ commutes with $\Pi_E (N)$, so that $N -\Pi_E (N)$ is normal.
Therefore, $N \in \cP_\eta \inc \cP_\beta \,$ and $T \in B_\Delta (\cP_\beta )$. \edem

\begin{teo}\label{es continua3}
The map $\Delta^\infty$ is continuous on $\matinv$.  
\end{teo}

\bdem  
By Remark \ref{reduccion}, it is enough to prove the continuity at the normal matrices of $\matinv$. 
Fix $N_0 \in \matinv $ a normal matrix. Let $\eps >0$, such that $B(N_0 \, , \eps)\inc \matinv$. 
We shall use the notatios of the previous statements relative to $N_0\,$. 
For $\beta  < \min \{\frac\eps2 , \eps_\mu \}$ small enough, we can extend to the 
compact set $\cP_\beta$ the distribution 
of subspaces $N \mapsto F_N$ and $\sub{\cE}{N}^s$ given by Corolary \ref{mitades} (for $\cP = \orbu{N_0}\,$), 
by using the functional calculus on the derivatives 
$T_M\Delta$, for $M \in \cP_\beta\,$ (see Remark \ref{CFD}).   
In this case, the subspaces $F_M$ and $\sub{\cE}{M}^s$ are $T_M\Delta$-invariant, 
$T_M\Delta \big|_{F_M} $ is near  $I_{F_M}$ and 
$\left\|\left. T_M\,\Delta\right|_{\sub{\cE}{M}^s}\right\|\leq \sub{k'}{N_0} <1$, 
for some $\sub{k}{N_0} < \sub{k'}{N_0} <1$. 

Observe 
that the set $\cP_\beta$ consists of fixed points for $\Delta$. 
Hence this distribution are in the hypothesis of  Proposition \ref{pelitosA}. 
Let $\rho>0 $ such that $\ds B(N_0 \, , \rho) \inc \mla{ \beta} \inc B_\Delta (\cP_\beta )\,$. 
Following the same steps of the proof of Proposition \ref{los w cortan}, but using 
Lemma  \ref{cdp2} instead of Lemma \ref{cdp}, we obtain that, if  $\rho$ 
is small enough then,  for every $T \in B(N_0 \, , \rho)$,   there exists 
$$
N_1 \in \Pi_E (\ewe^{ss}_T ) \cap \orbu{N_0} \cap B\big(N_0 \, , \frac\eps2 \, \big) \ .
$$
Let $S \in \ewe^{ss}_T $ such that $\Pi_E(S) = N_1\,$. Since $E(S) = E(N_1)$ is an orthogonal 
system of projectors, Proposition \ref{sumas} and Eq. \eqref{lasW} 
of Proposition \ref {pelitosA} assure that 
$$
\aluit{\infty} {T } = \aluit{\infty} {S } = N_2 
\peso {and} \Pi_E (N_2 ) = \Pi_E (S ) = N_1 \ .
$$
Since $T \in \mla{\beta}\,$, Lemma \ref{cosas} assures 
that 
$$\ds \|N_2 - N_1 \| = \|N_2 - \Pi_E(N_2) \| 
< \beta < \frac \eps2 \  .
$$ 
This shows  that $\aluit{\infty} { B(N_0 \, , \rho) \, } \inc  B(N_0 \, , \eps)$, 
i.e., that  $\Delta^{\infty}$ is continuous at $N_0\,$. 
\edem

\begin{rem} \label{conj}\rm 
By Remark \ref{fallo}, the Aluthge transform  fails to be  differentiable at 
every non invertible normal matrix. By this fact, we can not use the previous techniques 
for proving continuity of $\Delta^\infty$ on $\mat\setminus \matinv$. 
We conjecture that it is, indeed, continuous on $\mat$, but we have no proof 
for non invertible matrices. \EOE
\end{rem}

\section{Concluding remarks}
\subsection{Rate of convergence}
In \cite{[APS]} we proved that, if $T\in \mat$ is diagonalizable, then after some 
iterations the rate of  convergence of the sequence $\aluit{n}{T} $ becomes exponential. 
More precisely, 
for some $n_0\in\N$ and every $n\geq n_0 \,$, there exist $C>0$ and  
$0<\gamma <1$  such that $\|\aluit{n}{T}  - \aluit{\infty}{T}\| < C \gamma ^n $. 
This exponential rate depends on the spectrum of $T$. Actually, if $\la (T) = \la(D)$ for 
some diagonal matrix $D$, then  $\gamma = k_D\,$, the constant which appears in 
Theorem \ref{the key}. Using the formula for $k_D\,$, one can see that  
it  is closer to $1$ (so that the rate of convergence becomes slower) 
if the different eigenvalues are closer one to each other. 
 
These facts are not longer true if $T$ is not diagonalizable, since the rate of convergence
for such a $T$ depends on the rate of convergence for some 
$M \in  \ewe^{ss}_{ T} \cap \lentas $ (with the notation of Proposition \ref{los w cortan}), 
which can be much slower (and not exponential). Observe that the proof of the convergence 
of the sequence $\{\aluit{n}{M}\}$, given in Proposition \ref {en oded converge}, does not 
study the rate of convergence. It only shows that there exists an unique 
possible limit point for the sequence.

Nevertheless, using Proposition \ref{los w cortan} 
and Eq. \eqref{lasW}, 
it is easy to see that the system of projections  
$E(\aluit{n}{T}\,)$ converges to $E(\aluit{\infty}{T}\,)$ exponentially, because  
$E(M) = E(\aluit{\infty}{T}\,)$. 
As in the case of diagonalizable 
matrices the rate of convergence of the spectral projections depends on the 
spectrum of $T$, which agree with the spectrum of $M$.
Note that the spectrum of $T$ and the spectral projections of $M$ completely 
characterize the limit  $\aluit{\infty}{T}$. Indeed, 
if $\sigma(T) = \{ \mu_1, \dots , \mu_k\}$, then 
$$
\aluit{\infty}{T}=\aluit{\infty}{M} = \Pi_E(M) =  \sum_{j= 1}^k \mu_j \, E_j(M)
\ .
$$

\subsection{$\la$-Aluthge transform}
Given $\la \in (0,1)$ and a matrix $T\in\mat$ whose polar decomposition is $T=U|T|$, the $\la$-Aluthge transform of $T$ is defined by
$$
\alul{T}=|T|^\la U |T|^{1-\la}\ .
$$
All the results obtained in this paper are also true for the $\la$-Aluthge transform 
for every $\la \in (0,1)$, with almost the same proofs. 
Indeed, note that the basic results about Aluthge transform used throughout sections 
3 and 4 are Theorem \ref {the key} and those stated in subsection 2.1. 
All these results were extended to every $\la$-Aluthge transform   
(see \cite{[AMS]} and \cite {[APS2]}). The unique difference is that the constant
$k_{D}$ of Theorem \ref {the key} now depends on $\la$ (see Theorem 3.2.1 of \cite {[APS2]}). 
Anyway, the new constants are still lower than one for every $\la \in (0,1)$. Moreover, they 
are uniformly lower than one on compact subsets of $(0,1)$.

Another result which depends particularly on the  Aluthge transform is 
Proposition \ref{foto}, which is used to prove Corollary \ref{mitades}. 
Nevertheless, it is easy to see that both results 
are still true for every $\la \in (0,1)$.
On the other hand, the proof of Theorem \ref{es continua3} uses the same 
facts about the Aluthge transform. So that, it also remains true 
for $\Delta_\la\,$, for every $\la \in (0,1)$. 
We resume all these remarks in the following statement:
\begin{teo} \label{el de la} \rm 
For every $T \in \mat$ and $\la \in (0,1)$, the sequence $\alulit{n}{T} $ converges to a normal matrix 
$\alulit{\infty}{T} $. The map $T \mapsto \alulit{\infty}{T} $ is 
continuous on $\matinv$. \QED
\end{teo}

\noi We extend the conjecture given in Remark \ref{conj} to the following: 

\begin{conj} \rm 
The map $(0,1) \times \mat \ni (\la , T) \mapsto \alulit{\infty}{T}$ is continuous .
\end{conj}
Using the same ideas as in section 4 of \cite{[APS2]}, it can be proved that the above map is continuous if it 
is restricted to $(0,1) \times \matinv$. \EOE

\appendix{
\section{Appendix}}

\subsection {Proof of Proposition \ref{conitosA}}

Recall that $\cP $ is a compact set consisting of fixed points of $\Delta$ with two 
complementary, continuous and $T_N\Delta $ invariant distributions 
 $N \mapsto \sub{\cE}{N}^s $ and $N \mapsto \sub{\cF}{N}$ 
such that 
	\begin{equation}\label{menor que rho2}
	\left\|\left. T_N\,\Delta\right|_{\sub{\cE}{N}^s}\right\|< 1- \rho \peso{and} 
\left\|\left. (I - 	T_N \,\Delta\right)|_{\sub{\cF}{N}}\right\|<\frac \rho2 \ , \ \ N \in \cP \ , 
	\end{equation}
for some  $\rho  \in(0,1)$ which does not depend on $N$. 
The aim of the Proposition is to extend them
to distributions defined in some local basin of  $\cP$ with almost the same properties. 

The first step is to extend these distributions using 
the functional calculus: 
Fix $\eps>0$ such that  $\sigma(T_T\D ) \inc  B(1 , \, \frac \rho2 \, ) \cup B(0 ,  1-\rho) $, 
for 
every  
 $T\in B_\Delta (\cP)_\eps \,$.   
As in Remark \ref{CFD}, consider the spectral subspaces 
 $$
 F_T = R\big( \aleph _{B(1 , \, \frac \rho2 \,  )} (T_{T} \Delta  \,) \,\big)
 \peso{and} 
 \sub{E}{T} = R\big( \aleph _{B(0 , 1-\rho   )} (T_{T} \Delta  \,) \, )\ . 
$$
Observe that Eq. \eqref{menor que rho2}  assures that $E_N = \sub{\cE}{N}^s$ and $F_N = \sub{\cF}{N}$
for every $N \in \cP$. 
Since the functional calculus is smooth and $\cP$ is compact, we can assume that, for every 
$T \in B_\Delta (\cP)_\eps \,$, the angle between $F_T$ and $E_T$ is 
uniformly bounded from below, and $T_T\D$ satisfies inequalities as in Eq. \eqref{menor que rho2}, when 
it is restricted to $E_T$ and $F_T\,$.  
Let us take the 
cones $C_T=C(\alpha, E_{T})$  of size  $\alpha$ in the direction $E_{T}\,$.  
For every small $\alpha$, we can assume that 
\ben
\item [a)] There exists $\gamma>0$ such  that $C_T \cap C(\gamma , F_T) = \{0\}$ 
for every $T\in B_\Delta (\cP)_\eps \,$. 
\item [b)] Every subspace $E'_{T} \inc C_T $ with $\dim E'_{T} = \dim E_{T}\,$ satisfies 
inequalities as in Eq \eqref{menor que rho2}. 
\een

\begin{claim} \label{Cl1}
There exist positive constants $\la_0<1$ and $\al>0$ such that,  
 if $\eps$ is a small enough, then for every   $T \in B_\Delta (\cP)_\eps $ it holds that 
  $[T_{T}\D]^{-1}(C_{\D(T)})$ is a cone of size not greater than $\la_0 \, \alpha$ 
 inside $C_{T}\,$.
\end{claim}

\noindent {\it Proof of the claim:} First observe that, by 
the properties of the subspaces $E_T$ and $F_T\,$, there exist
$\la_1<1$ and $\al>0$ such that,  for every   $T \in B_\Delta (\cP)_\eps $ it holds that 
  $[T_{T}\D]^{-1}(C_{T}) \inc  C(\la _1 \, \alpha, E_{T})$ which 
  is a cone of size $\la_1 \, \alpha$  inside $C_{T}\,$.

Take $T\in B_\Delta (\cP)_\eps $ and its image $\D(T) \in B_\Delta (\cP)_\eps \,$. 
Observe that  $\D$ commutes with unitary conjugations, and 
$\D$ is uniformly continuous on compact sets. Hence, if 
$\eps$ is taken small enough, then 
 $T$ is arbitrarily (and uniformly) close to $\D(T)$ for every $T \in B_\Delta(\cP)_\eps \,$. 
 Therefore, $E_{T} $  is arbitrarily and uniformly 
close to $E_{\D(T)}\,$, and the same occurs between $C_{\D(T)}$ and $C_{T}\,$. 
Putting all together, it follows that 
 $$
 [T_T \D]^{-1}(C_{\D(T)})
 \sim  
 [T_{T}\D]^{-1}(C_{T}) \inc  C(\la _1 \, \alpha, E_{T}) \ .
 $$ 
 Therefore, there exists $\la_1<\la_0 <1$ such that  
$[T\D]^{-1}(C_{\D(T)})$ is a cone of size not greater than $\la_0 \, \alpha$  inside  
 $C_T\,$. This completes the proof of the Claim. 

\medskip
\noi
It is easy to see that the Claim  implies that, if $C$ is a cone of size $\beta <\alpha$
inside $C_{\D(T)}$ and of the same dimension, then 
$[T_{\D(T)}\D]^{-1}(C)$ is a cone of size not greater than $\la_0 \beta$ inside $C_{T}\,$.  
For each $T\in B_\Delta (\cP)_\eps\,$, consider the sequence  
$\{\D^{n}(T)\}_{n\in\N}$ and the sequence of cones 
$$
C_1 = [T_{T}\D]^{-1}(C_{\D(T)}) 
\peso{and} 
\{C_n \}_{n\in\N} = \{[T_{T}\D^n]^{-1}(C_{\D^n(T)})\}_{n\in\N} \ 
$$ 
in $T_T\mat$.  The following facts hold: 
For every $n \in \N$, 
$$
\barr{rl}
C_{n+1} &= [T_{T}\D^{n}]^{-1}\big( \, [T_{\D^{n} (T)}\D]^{-1}  C_{\D^{n+1}(T)}\, \big)
\\ & \inc [T_{T}\D^{n}]^{-1} (C_{\D^{n}(T)}) = C_{n}\ .
\earr
$$
Therefore $C_{n+1} \inc C_n \inc C_1 \inc C_T \,$ and  
every $C_n$ is a cone of size not greater than  $\la_0 ^n \, \al\,$.   
An easy  argument of dimensions shows that every set $C_n$ contains a subspace of 
dimension equal to $\dim E_{T}\,$ (even if the derivarives $T_{\D^n (T)} \D$ are not bijective). 
Therefore, 
$$
\cE^s_T \ := \ \bigcap_{n \in \N}  C_n =\bigcap_{n \in \N} [T_{\D^n(T)}\D]^{-n}(C_{\D^n(T)})
$$ 
is a a well defined unique direction, and $\dim  \cE^s_{T} = \dim E_{T}\,$. 
Observe that the direction is invariant and  $\cE^s_{T} \inc C_T\,$, 
and so it is contracted by $T\D$.
Take $\cF_T = F_T\,$, $T \in B_\Delta (\cP)_\eps\,$, which is continuous by construction. 
The  continuity of $\cE^s_{T}$  follows from the fact that this subbundle is invariant and 
uniformly contracted for any forward iterate and from the uniqueness of a 
subbundle (with maximal dimension) exhibiting  these properties.  Finally, the subspaces 
$\cE^s_{T}$ and $\cF_T$ satisfy Eq \eqref{menor que rho2} by construction. 
\edem

\subsection{Proof of Proposition \ref{los w cortan}}

In this section we shall use the following notations: If $D\in \mat$ is an invertible 
diagonal matrix then $\cP = \orbu{D}$ and, for every $T \in  B_\Delta (\cP)_\eps\,$,  by means of
$\ewe^{ss} _{ T} = \ewe^{ss} ( T)   : (-1, 1)^m \to B_\Delta (\cP)$ we denote the maps 
given by Proposition \ref{pelitosA}. The invariant manifolds will be denoted by 
$\ewe^{ss} _{ T} \big(\, (-1, 1)^m \, \big)$. Finally, for a sake of simplicity, for every
$t>0$, $\cubo{t}$ denotes the $m$-dimensional cube $(-t, t)^m$.

Observe that $\ewe^{ss}_{ T}(\cubo{1}) $ intersects $\lentas$ if and only if 
$\Pi_E(\ewe^{ss}_{T}(\cubo{1})\,)$ intersects $\cP$. The proof of Proposition \ref{los w cortan} uses this remark and it is based on some well known results about transversal intersections, 
using that $\Pi_E (\ewe^{ss}_{ T}(\cubo{1})\,)$ is ``$C^2$-close" to another 
mainfold ($\ewe^{ss}_{ N}(\cubo{1})$ for some $N\in \cP$ near $T$)
which intersecs transversally $\cP$, both contained in $\orb{D}$. We give a proof adapted to our case, divided into three lemmas: We begin with the following classical result (see for example \cite[pg.36]{[H]}).

\begin{lem}\label{Hirsch} \rm
Let $U\subseteq \R^m$ be an open set and $W\subseteq U$ an open set with compact 
closure $\overline{W}\subseteq U$. Let $M \inc \R^n$ be a smooth submanifold and  
$f:U\to M$ a $C^1$ embedding. There exists $\eps>0$ such that, if 
$$
g:U\to M \ \mbox{is $C^1$ , \ } 
\|T_xg-T_xf\|<\eps \peso{and} \|g(x)-f(x)\|<\eps
$$
for every $x\in W$, then $\left. g\right|_{W}$ is an embedding. \QED
\end{lem}

\begin{lem}\label{cdp}
Let $D$ and $\cP$ be as in Proposition \ref{los w cortan}.  
Then there is $\eta<\eps$ such that the map  \rm
$$
\cV: B_\Delta (\cP)_\eta   \to \mbox{Emb}^2\big( \, \cubo{\frac{1}{2}}\, ,\orb{D} \, \big) 
\mbox{ given by } 
\cV_T=\left.\Pi_E\circ\ewe^{ss}_T\right|_{\cubo{\frac{1}{2}}}  \ ,
$$ \it
is well defined and continuous with respect to the the $C^2$ 
topology of \rm $\mbox{Emb}^2\big( \, \cubo{\frac{1}{2}}\, ,\orb{D} \, \big)$.
\end{lem}

\bdem
Consider the map  $\widetilde\cV: B_\Delta (\cP)_\eps   \to 
C^2(\cubo{1}\, ,\orb{D}\, )$ given by $\widetilde\cV_T=\Pi_E\circ\ewe^{ss}_T\,$. 
By Proposition \ref{pelitosA} and Remark \ref{E y PIE}, $\widetilde\cV$ is  
well defined and continuous, 
 if $C^2(\cubo{1}\, ,\orb{D}\, )$ 
is endowed with the $C^2$ topology. 
Observe that $\ewe^{ss}_N $ takes values in $\orb{D}$ for every $N\in\cP$. Indeed, 
this follows by Cor. 3.1.2 of \cite{[APS]}, or by rewriting the proof 
of Proposition \ref {pelitosA} inside $\orb{D}$ in this case. 
Therefore, 
$\widetilde\cV_N=\ewe^{ss}_N$ for every $N\in\cP$, because 
$\Pi_E$ is the identity on $\orb{D}$. 

Given $N\in \cP$, Lemma \ref{Hirsch} assures that  there exists $\eps_N$ such that, 
if $\ete: \cubo{1}\to\mat$ is a $C^1$ map which satisfies that 
\begin{equation}\label{seras lo que debas ser}
\|T_x\ewe^s_N-T_x\ete\|<\eps_N \peso{and} \|\ewe^s_N(x)-\ete(x)\|<\eps_N \ ,
\end{equation} 
for every $x\in \cubo{\frac12}\,$, 
then $\ete\big|_{\cubo{\frac12}}$ is an embedding. By the continuity of 
$\widetilde\cV$, there is a neighborhood $\cU_N$ of $N$ in $B_\Delta (\cP)_\eps$ 
such that, for every $T\in\cU_N\,$, the map $\widetilde\cV_T$ satisfies 
\eqref{seras lo que debas ser}. Take $\eta>0$ such that 
$B_\Delta (\cP)_\eta\subseteq \bigcup_{N\in\cP}\cU_{N}\,$. 
Then, 
$\cV(T) = \widetilde\cV _T \big|_{\cubo{\frac12}} \in \mbox{Emb}^2(\cubo{\frac12}\, ,\orb{D}\, )$, 
for every $T\in B_\Delta (\cP)_\eta\,$ i.e., $\cV$ is well defined. The continuity of $\cV$ follows from the 
fact that both $\widetilde\cV$ and the restriction map 
$\ete \mapsto \ete\big|_{\cubo{\frac12}}$ are continuous with respect to the the $C^2$ 
topology. 
\edem

\begin{lem}\label{interseccion transversal}
Let $D$ and $\cP$ be as in Proposition \ref{los w cortan}.  
Given $N_0\in \cP$ and $\eps>0$, 
there exists a $C^2$-neighborhood $\Omega$ of $\ewe^{ss}_{N_0}\big|_{\cubo{\frac12}}$ 
in the space $\mbox{\rm Emb}^2(\cubo{\frac12} ,\orb{D}\,)$,  
such that  $\ete (\cubo{\frac12})$  intersects the submanifold $\cP$ at a 
point $N\in B(N_0,\eps)$, for every $\ete\in \Omega$.
\end{lem}
\bdem
Let $(\cU_{N_0}\, , \fii)$ be a chart in $\orb{D}$ such that 
$N_0\in \cU_{N_0}\subseteq B(\eps, N_0)$, 
$\fii(N_0)=0$, and $\fii(\cP\cap \cU_{N_0})=\fii(\cU_{N_0})\cap (\{0\}\oplus\R^{n-m})$, 
where $\R^n\simeq\R^m\oplus\R^{n-m}$. Let $P$ denote the orthogonal projection from $\R^n$ 
onto $\R^m\oplus \{0\}$. 

By Proposition \ref{pelitosA},  
the intersection $\ewe^{ss}_{N_0}(\cubo{\frac12})\cap\orbu{D}=\{N_0\}$ is transversal.  
Then, there exist $\delta\in (0,1/2)$ and a $C^2$-neighborhood $\Omega_0$ of 
$\left.\ewe^{ss}_{N_0}\right|_{\cubo{\frac12}}$ in 
$\mbox{Emb}^2(\cubo{\frac12}\, , \orb{D}\,)$ such that, 
for every $\ete\in\Omega_0\,$, 
\begin{enumerate}
	\item $\ete(\,\overline{\cubo{\delta}}\,)\subseteq \cU_{N_0}$;
  \item  $\ker  P  \oplus T_M \widetilde\ete(\cubo{\delta}) =\R^n$, where 
 $\widetilde\ete=\left.\fii\circ\ete\right|_{\cubo{\delta}}$ and 
 $M \in \widetilde\ete(\cubo{\delta})$. 
 \item The angle between $\ker  P $ and $ T_M \widetilde\ete(\cubo{\delta})$ is uniformly bounded from below.
\end{enumerate}
Note that items (2) and  (3) imply that,  
$\R^k\oplus\{0\}$. On the other hand, item (2) also implies that 
for every $\ete\in\Omega_0$ and every $M \in  \widetilde\ete(\cubo{\delta})\,
$,  
the linear map $P$ acting on $ T_M \widetilde\ete(\cubo{\delta})$ is uniformly bounded from below. 
On the other hand, the norm of the second derivative of 
$P\circ \widetilde\ete$ is bounded on $\cubo{\frac{\delta}{2}} $. 
Hence there exists $\mu >0$ so that,  
for every $ M \in  \widetilde\ete  \left(\cubo{\frac{\delta}{2}}\right)$, 
\begin{equation}\label{little balls}
 B(P(M), \mu ) \inc P \big(\widetilde\ete(\cubo{\delta}) \, \big)  \ .
\end{equation}
Take $\Omega\subseteq \Omega_0$ such that 
$\|\widetilde{\ewe}^{ss}_{N_0}(x)-\widetilde\ete(x)\|<\mu/2$ for every $\ete\in \Omega$ 
and every $x\in \cubo{\delta}\, $, where $\widetilde{\ewe}^{ss}_{N_0}=\fii\circ \ewe^{ss}_{N_0}\big|_{\cubo{\delta}}\,$.  As $\widetilde{\cW}^{ss}_{N_0}(0)=0$, 
Eq. \eqref{little balls}
implies that  $0 \in  P \big(\widetilde\ete(\cubo{\delta}) \, \big)$,  for every  $\ete\in \Omega$. 

Thus $\ete\cap \orbu{D}\cap \cU_{N_0}\neq\varnothing$, because 
$\ete(\, \overline{\cubo{\delta}} \, )\subseteq \cU_{N_0}\,$. In particular, $\ete(\cubo{\frac12})$ intersects 
the submanifold $\cP$ transversally at a point $N\in \cU_{N_0} \inc B(N_0,\eps)$. 
\edem

\medskip
\noi \it Proof of Proposition \ref{los w cortan}: \rm 
Given $N\in\cP$, Lemma \ref{interseccion transversal} assures that 
there is a $C^2$-neighborhood $\Omega_N$ of 
$\ewe^{ss}_N\big|_{\cubo{\frac12}}$ in $\mbox{Emb}^2(\cubo{\frac12}\, ,\orb{D})$ such that 
$\ete(\cubo{\frac12}) \cap \cP\neq \varnothing$ for every $\ete\in \Omega_N\,$. Let $\cV$ be the 
function defined in Lemma \ref{cdp}, and let $\cU_N=\cV^{-1}(\Omega_N)\,$. 
Define $\cU_\cP=\bigcup_{N\in\cP}\cU_N\,$. Therefore, $\cU_\cP$ is an 
open neighborhood of $\cP$ contained in $B_{\Delta}(\cP)\,$. Since $\cP$ is compact, 
there exists $0<\eta <\eps$ such that $B_\Delta (\cP)_\eta \inc \cU_\cP\, $. 
Then, for every 
$T \in B_\Delta (\cP)_\eta\,$,  $\Pi_E \big( \ewe^{ss}_T (\cubo{\frac12})\, \big)$ itersects $\cP$.  
By Proposition \ref {pelitosA}, $ \ewe^{ss}_T (\cubo{\frac12}) \inc B_{\Delta}(\cP)\,$.
So that $ \ewe^{ss}_T (\cubo{\frac12}) \cap  \lentas \neq \varnothing\,$. 
\QED

\medskip
\noi 
The proof of the next result, which is 
used in the proof of the continuity of the limit function $\Delta^\infty$,  
follows the same lines as the proof of  Lemma \ref{cdp}.

\begin{lem}\label{cdp2} \rm 
Let $N_0\in\mat$ be a normal matrix, $\cP_\beta$ as in Definition \ref{pebeta} 
and $\ewe^{ss}: B_\Delta (\cP_\beta)_\eps \to \mbox{\rm Emb}^2(\cubo{1}\, ,B_\Delta (\cP) )$  
the prelamination given by Proposition \ref{pelitosA}. If $\Pi_E$ is defined with 
respect to the spectrum of $N_0\,$, then there exists $\eta<\beta$ so that the map 
$$
\cV: B_\Delta (\cP_\eta)_\eta   \to \mbox{Emb}^2\big( \, \cubo{\frac12}\,,\orb{N_0} \, \big)
$$
given by $\cV_T= \Pi_E\circ\ewe^{ss}_T\big|_{\cubo{\frac12}} 
$ 
is well defined and continuous with the $C^2$ 
topology of \rm $\mbox{Emb}^2\big( \, \cubo{\frac12},\orb{N_0} \, \big)$.  \QED
\end{lem}

\bigskip
\fontsize {9}{11}\selectfont

\noi {\bf Jorge Antezana and Demetrio Stojanoff}

\noi Depto. de Matem\'atica, FCE-UNLP,  La Plata, Argentina
and IAM-CONICET  

\noi e-mail: antezana@mate.unlp.edu.ar and demetrio@mate.unlp.edu.ar

\medskip

\noi {\bf Enrique R. Pujals}

\noi Instituto Nacional de Matem\'atica Pura y Aplicada (IMPA), Rio de Janeiro, Brasil.

\noi e-mail: enrique@impa.br

\end{document}